\documentclass[12pt]{amsart}

\makeatletter
\@namedef{subjclassname@2020}{\textup{2020} Mathematics Subject Classification}
\makeatother
\newcommand{\er}{\eqref}

\newcommand{\sodd}{\wh{\sig_{}}}

\newcommand{\tanz}{\tan z}
\newcommand{\cotz}{\cot z}
\newcommand{\piit}{\pii \tau}

\newcommand{\mg}{\infty}

\newcommand{\pii}{\pi i}

\newcommand{\pz}{\pi z}
\newcommand{\piz}{\pi z}

\renewcommand{\r}{\mathbf{r}}

\newcommand{\R}{\mathcal{R}}

\newcommand{\dsum}{\di\sum}

\newcommand{\wh}{\ensuremath{\widehat}}
\newcommand{\el}{\ensuremath{\ell}}

\DeclareFontFamily{U}{mathx}{\hyphenchar\font45}
\DeclareFontShape{U}{mathx}{m}{n}{
      <5> <6> <7> <8> <9> <10>
      <10.95> <12> <14.4> <17.28> <20.74> <24.88>
      mathx10
      }{}
\DeclareSymbolFont{mathx}{U}{mathx}{m}{n}
\DeclareFontSubstitution{U}{mathx}{m}{n}
\DeclareMathAccent{\widecheck}{0}{mathx}{"71}

\renewcommand{\kill}[1]{}
\newcommand{\dummy}[1]{\mbox{}}

\makeatletter
\newcommand{\xequal}[2][]{\ext@arrow 0055{\equalfill@}{#1}{#2}}
\def\equalfill@{\arrowfill@\Relbar\Relbar\Relbar}
\makeatother

\newcommand{\mto}{\mapsto}

\newcommand{\Set}[2]{\ensuremath{\left\{{#1}\,\middle|\,{#2}\right\}}}

\renewcommand{\k}{\ensuremath{\ol{\mathrm{P}}}}

\newcommand{\hou}[3]{{#1}\equiv {#2}\pmod{#3}}

\newcommand{\h}{\hline}

\renewcommand{\k}[1]{\ensuremath{\left({#1}\right)}}

\newcommand{\ds}{\dots}

\newcommand{\bca}{\begin{cases}}
\newcommand{\eca}{\end{cases}}

\renewcommand{\th}{\ensuremath{\theta}}

\newcommand{\mug}{\ensuremath{\infty}}
\newcommand{\der}[1]{\ensuremath{\lef({#1}\ri)'}}

\newcommand{\dprod}{\di\prod}

\newcommand{\del}{\ensuremath{\delta}}

\newcommand{\Gam}{\Gamma}

\newcommand{\ff}[2]{\ensuremath{\di\fr{#1}{#2}}}

\newcommand{\s}[1]{\ensuremath{\di\int{#1}\,dx}}

\newcommand{\bpic}{\begin{picture}}\newcommand{\epic}{\end{picture}}

\newcommand{\beda}{\begin{edaenumerate}}
\newcommand{\eeda}{\end{edaenumerate}}

\newcommand{\cd}{\cdots}

\newcommand{\sh}[1]{\shadowbox{#1}}

\newcommand{\q}{\quad}

\newcommand{\ta}[1]{\ensuremath{\tan\left({#1}\right)}}

\newcommand{\bq}{\begin{quote}}\newcommand{\eq}{\end{quote}}

\newcommand{\sig}{\sigma}

\newcommand{\ti}{\times}

\newcommand{\be}{\begin{enumerate}}\newcommand{\ee}{\end{enumerate}}
\newcommand{\bce}{\begin{center}}\newcommand{\ece}{\end{center}}
\newcommand{\bde}{\begin{description}}\newcommand{\ede}{\end{description}}
\newcommand{\bri}{\begin{flushright}}\newcommand{\eri}{\end{flushright}}
\newcommand{\bb}{\begin{block}}\newcommand{\eb}{\end{block}}
\newcommand{\bt}{\begin{thm}}\newcommand{\et}{\end{thm}}
\newcommand{\bpf}{\begin{proof}}\newcommand{\epf}{\end{proof}}
\newcommand{\bex}{\begin{ex}}\newcommand{\eex}{\end{ex}}
\newcommand{\bexr}{\begin{exr}}\newcommand{\eexr}{\end{exr}}
\newcommand{\bft}{\begin{fact}}\newcommand{\eft}{\end{fact}}
\newcommand{\brk}{\begin{rmk}}\newcommand{\erk}{\end{rmk}}
\newcommand{\ba}{\begin{align*}}\newcommand{\ea}{\end{align*}}
\newcommand{\bexe}{\begin{exe}}\newcommand{\eexe}{\end{exe}}

\newcommand{\bit}{\begin{itemize}}\newcommand{\eit}{\end{itemize}}

\newcommand{\bcm}{}

\newcommand{\hf}{\hfill}
\newcommand{\fr}{\frac}
\newcommand{\nn}{\ensuremath{\mathbb{N}}}

\newcommand{\zz}{\ensuremath{\mathbb{Z}}}

\newcommand{\bd}{\begin{defn}}\newcommand{\ed}{\end{defn}}
\newcommand{\bp}{\begin{prop}}\newcommand{\ep}{\end{prop}}
\newcommand{\p}{\ensuremath{\pi}}
\newcommand{\eh}{\emph}\newcommand{\al}{\alpha}

\newcommand{\mb}{\mbox}
\newcommand{\te}{\text}
\newcommand{\wt}{\widetilde}
\newcommand{\lef}{\left}\newcommand{\ri}{\right}
\renewcommand{\l}{\left}\renewcommand{\r}{\right}

\newcommand{\di}{\displaystyle}

\renewcommand{\d}{\ensuremath{\bm{d}}}

\newcommand{\z}{\ensuremath{\bm{z}}}
\newcommand{\np}{\newpage}

\renewcommand{\s}{\sigma}

\renewcommand{\d}{\delta}
\renewcommand{\P}{\mathcal{P}}
\newcommand{\Q}{\mathcal{Q}}

\renewcommand{\d}{\delta}

\renewcommand{\ta}{\tau}

\usepackage[dvipdfmx]{graphicx}
\usepackage[dvipsnames]{xcolor}
\usepackage{float,afterpage}
\usepackage{asymptote,layout}
\usepackage{wrapfig,epic}
\usepackage{amsmath,amssymb,cases,pict2e}

\usepackage{multicol}
\usepackage{tikz}

\usepackage{geometry}
\usepackage{exscale,latexsym,bm}
\usepackage{amssymb,enumerate,amsmath,amsthm,amsfonts}
\usepackage{verbatim,fancybox,wasysym,fancyhdr,type1cm}
\usepackage[frame,all,poly,curve,knot,arrow]{xy}
\usepackage{colortbl}
\usepackage{boxedminipage}
\usepackage{multirow}

\graphicspath{{./figs/}}
\theoremstyle{definition}
\newtheorem{thm}{Theorem}[section]

\newtheorem{lem}[thm]{Lemma}
\newtheorem{prop}[thm]{Proposition}\newtheorem{cor}[thm]{Corollary}
\newtheorem{cj}[thm]{Conjecture}

\newtheorem{exr}[thm]{Exercise}

\newtheorem{ex}[thm]{Example}

\newtheorem{defn}[thm]{Definition}\newtheorem{rmk}[thm]{Remark}
\newtheorem{fact}[thm]{Fact}
\newtheorem{block}[thm]{}
\newtheorem*{exe}{Exercise}

\renewcommand{\l}{\zeta^{*}}
\renewcommand{\h}{\hline}
\renewcommand{\arraystretch}{1.5}

\renewcommand{\P}{\mathbf{P}}
\renewcommand{\z}{\zeta}

\newcommand{\tf}{\tfrac}
\newcommand{\tff}{\tfrac}
\newcommand{\innn}{\in\mathbb{N}}

\newcommand{\sinz}{\sin z}
\newcommand{\cosz}{\cos z}
\newcommand{\qqm}{(q;q)_{\mug}}
\renewcommand{\sh}{\widehat{\sigma}}

\renewcommand{\P}{\mathcal{P}}
\newcommand{\E}{\mathcal{E}}
\newcommand{\incc}{\in\cc}
\renewcommand{\del}{\partial^{*}}

\newcommand{\cc}{\ensuremath{\mathbb{C}}}
\newcommand{\qq}{\ensuremath{\mathbb{Q}}}

\begin{document}

\title{Ramanujan-Shen's differential equations 
for Eisenstein series of level 2
}
\author{Masato Kobayashi}
\date{\today}

\address{Masato Kobayashi\\
Department of Engineering\\
Kanagawa University, 3-27-1 Rokkaku-bashi, Yokohama 221-8686, Japan.}
\email{masato210@gmail.com}

\maketitle
\begin{abstract}
Ramanujan (1916) and Shen (1999) discovered differential equations for classical Eisenstein series. 
Motivated by them, we derive new differential equations for Eisenstein series of level 2 from the second kind of Jacobi theta function. This gives a new characterization of a system of differential equations by Ablowitz-Chakravarty-Hahn (2006), Hahn (2008), Kaneko-Koike (2003), Maier (2011), Nidelan (2022) and Toh (2011). As application, we show some arithmetic results on Ramanujan's tau function. 
%
\end{abstract}

\renewcommand{\ta}{\tau}

\renewcommand{\P}{\mathcal{P}}
\renewcommand{\Q}{\mathcal{Q}}
\renewcommand{\R}{\mathcal{R}}
\renewcommand{\E}{\mathcal{E}}

\section{Introduction}
\let\thefootnote\relax\footnote{This is 
M. Kobayashi, Ramanujan-Shen's differential equations for Eisenstein  series of level 2. Res. number theory 10, 41 (2024).}

\subsection{Ramanujan's paper in 1916}

Our motivation of this research is 
Ramanujan's work \cite{ramanujan} on the classical \eh{Eisenstein series}. For $n, s\innn$, let 
$\s_{s}(n)=\dsum_{d|n}d^{s}$ and 
for $q\incc$, $|q|<1$, set 
\begin{align*}
	P(q)&=1-24
	\dsum_{n=1}^{\mug}\s_{1}(n)q^{n},
	\\Q(q)&=1+240
	\dsum_{n=1}^{\mug}\s_{3}(n)q^{n},
	\\R(q)&=1-504
	\dsum_{n=1}^{\mug}\s_{5}(n)q^{n}.
\end{align*}
\begin{thm}[Ramanujan's differential equations]
The following holds.
\label{t1}
\begin{eqnarray*}
	qP'&=\ff{P^{2}-Q}{12}.
	\\qQ'&=\ff{PQ-R}{3}.
	\\qR'&=\ff{PR-Q^{2}}{2}.
\end{eqnarray*}
\end{thm}
Ramanujan proceeds as follows. 
For $x\in\cc, |q|<1$, set 
\begin{align*}
	{S}&=\ff{1}{4}\cot \ff{x}{2}+
	\dsum_{n=1}^{\mug}\ff{q^{n}}{1-q^{n}}\sin nx,
	\\{T}_{1}&=
	\k{\ff{1}{4}\cot \ff{x}{2}}^{2}
		+\dsum_{n=1}^{\mug}\ff{q^{n}}{(1-q^{n})^{2}}\cos nx,
	\\{T}_{2}&=
	\ff{1}{2}
	\dsum_{n=1}^{\mug}\ff{nq^{n}}{1-q^{n}}
	(1-\cos nx).
\end{align*}
Then he claimed $S^{2}=T_{1}+T_{2}$. 
Further, he derived 
$qP'=\tff{P^{2}-Q}{12}$ from coefficients of 
$x^{2}$ in that identity. 
Shen \cite{shen} extended this idea with \eh{Jacobi theta functions} as we will go into very details in Section 2.

\subsection{Main results}
Our main result is Theorem \ref{t3} which 
proves an infinite sequence of differential equations for Eisenstein series $(E_{2k}^{*})$ of level 2. This is an analog of Theorem \ref{t1} and \ref{t2}. 
In Sections 3 and 4, we will present its consequences 
(Theorems \ref{t314}, \ref{t10}, Corollaries \ref{c2}, \ref{c10}) on \eh{Ramanujan tau function}. 
Also, we will show a new closed system of differential equations (Proposition \ref{p4})
\[
\del A=-\ff{1}{4}(A^{2}+B), \q
\del B=-BC, \q
\del C=-\ff{1}{2}B 
\]
with 
$A=E_{2}^{*}, B=E_{4}^{*}, 
C=\tf{E_{6}^{*}}{E_{4}^{*}}$ 
and $\del$ certain differential operator as we explain later. 
See Ablowitz-Chakravarty-Hahn (2006) \cite{ach}, 
Hahn (2008) \cite{hahn2}, Kaneko-Koike (2003) \cite{kk},  
Maier (2011) \cite{maier}, 
Nidelan (2022) \cite{nikdelan}
 and Toh (2011) \cite{toh} for other equivalent descriptions of such a system.

\section{Ramanujan-Shen's differential equations}

\subsection{Preliminaries on Riemann zeta}

Define signed Bernoulli numbers $(B_{n})_{n\ge0}$ by 
\[
\ff{x}{e^{x}-1}=
\dsum_{n=0}^{\mug}\ff{B_{n}}{n!}x^{n}.
\]
The first few values of $B_{n}$ are 
\[
B_{0}=1, B_{1}=-\ff{1}{2}, 
B_{2}=\ff{1}{6}, B_{4}=-\ff{1}{30}, 
B_{6}=\ff{1}{42},\]
\[B_{8}=-\ff{1}{30}, B_{10}=\ff{5}{66}, 
B_{12}=-\ff{691}{2730}, 
B_{14}=\ff{7}{6}
\]
and $B_{2k+1}=0$ $(k\ge1)$. Euler proved 
\begin{eqnarray*}\label{z2k}
\z(2k)=-\ff{1}{2}\ff{(2\pii)^{2k}}{(2k)!}B_{2k}
\q (k\ge0).
\end{eqnarray*}

\begin{fact}\label{f1}
Both of the following hold.
\begin{enumerate}
\item For $|z|<1$, 
\begin{eqnarray}\label{zcot}
\pi z\cot \pi z=
-2\dsum_{k=0}^{\mug}
\z(2k)z^{2k}.
\end{eqnarray}
\item For $m\ge2$, 
\begin{eqnarray*}\label{z2m}
\z(2m)=\ff{2}{2m+1}
\dsum_{k=1}^{m-1}\z(2k)\z(2m-2k).
\end{eqnarray*}
\end{enumerate}
\end{fact}
These results follow from 
Maclaurin series of $z\cot z$ and 
the differential equation
$(\cotz)'=-1-\cot^{2}z$. Since proofs are elementary, we omit such details here. Note that although we deal with various infinite sums below, sometimes we will not worry about those radii of convergence if it is elementary.
Now let us introduce an odd analog of $\z(2k)$.
\begin{defn}
For $k\ge0$, define 
$\l(2k)=(1-2^{-2k})\z(2k)$.
\end{defn}
%


\begin{prop}We have the following.
\begin{enumerate}
\item For $|z|<1$, 
\begin{eqnarray}\label{tan}
-\ff{\pi z}{2}\tan
{\ff{\pz}{2}}=
-2\dsum_{k=0}^{\mug}
\l(2k)z^{2k}.
\end{eqnarray}
\item For $m\ge2$, 
\begin{eqnarray}
\label{z*2m}
\l(2m)=\ff{2}{2m-1} \sum_{k=1}^{m-1}\l(2k)\l(2m-2k).
\end{eqnarray}
\end{enumerate}
\end{prop}
\begin{proof}
Use Maclaurin series of $z\tan z$ and 
$(\tanz)'=1+\tan^{2}z$.
\end{proof}

\subsection{Shen's differential equation}

Shen \cite{shen} explained Ramanujan's results 
more abstractly in terms of the following differential equation.
\begin{eqnarray}
\label{sdf}
\k{\ff{f'}{f}(z)}^{2}=
\ff{f''}{f}(z)-
\k{\ff{f'}{f}(z)}'.
\end{eqnarray}
Here let us call this 
\eh{Shen's differential equation}. 
Note that 
$(\cot z)^{2}=-1-(\cot z)'$ shows 
$f(z)=\sinz$ is a solution of 
\er{sdf} and 
so is $\cosz=\sin\k{z+\tff{\p}{2}}$.
Shen brought this idea into \eh{Jacobi theta functions}.
For $\ta\in \mathbb{H}=\{\ta\in \cc\mid \text{Im}(\ta)>0\}$, 
we always assume $q=e^{2\pii \ta}$. 

\begin{defn}
For $z\incc, \ta\in\mathbb{H}$, define 
the first kind of Jacobi theta function by 
\[
\th_{1}(z|\ta)=
2
\dsum_{n=0}^{\mug}(-1)^{n}q^{(2n+1)^{2}/8}\sin(2n+1)z.
\]
\end{defn}

\begin{fact}[{\cite[p.301-302]{shen}}]Both of the following hold.
\begin{enumerate}
	\item The zeros of $\th_{1}$ are 
$m\p+n\ta\p, m, n\in \zz$; all of them are simple.
	\item For each $\tau\in \mathbb{H}$, 
the function $f(z)=\th_{1}(z|\ta)$ satisfies \eqref{sdf}.
\end{enumerate}

\end{fact}

\subsection{Ramanujan-Shen's differential equations}

Let $\s_{s}(n)=\sum_{d|n}d^{s}$ as before. Often, write $\s_{1}(n)=\s(n)$.
\begin{defn}
For $k\ge1$, set 
\[
G_{2k}(\tau)=
2\z(2k)+
2\ff{(2\p i)^{2k}}{(2k-1)!}
\dsum_{n=1}^{\mug} \s_{2k-1}(n)q^{n}
\]
and $G_{0}(\ta)=-1$ for all $\ta\in \mathbb{H}$ for convenience.
\end{defn}

\begin{fact}[{\cite[p.305 (20)]{shen}}]
We have 
\begin{eqnarray}\label{ldth1}
\p\ff{\th_{1}'}{\th_{1}}(\piz|\ta)=-
\dsum_{k=0}^{\mug}G_{2k}(\tau)z^{2k-1}.
\end{eqnarray}
\end{fact}
We now see striking similarity between this and \er{zcot}.
%

The point of Shen's discussion is that 
Ramanujan's identity ${S}^{2}={T}_{1}+{T}_{2}$
is equivalent to that 
$f(z)=\th_{1}(\piz|\ta)\, (x=2z)$ satisfies 
\eqref{sdf}.

\begin{defn}
For $k\ge0$, set 
\[
E_{2k}(\tau)=\ff{G_{2k}(\ta)}{2\z(2k)}
\]
and $\s_{2k-1}(0)=(-\tff{4k}{B_{2k}})^{-1}$.
\end{defn}
\begin{rmk}
Hereafter, whenever there is no confusion, we just write 
$E_{2k}$ (or even $E_{2k}(q)$ which is abuse of notation, though) to mean $E_{2k}(\tau)$. 
Usually we do not talk about $E_{0}$. 
However, introducing $E_{0}=1$ is helpful for the argument on Hankel determinants of $(E_{2k})$ later on.
\end{rmk}

For example, $E_{2}=P, E_{4}=Q, E_{6}=R$ 
and 
\[
E_{8}=
1+{480}
\dsum_{n=1}^{\mug}\s_{7}(n)q^{n}, \q 
E_{10}=
1-264
\dsum_{n=1}^{\mug}\s_{9}(n)q^{n}.
\]
\begin{fact}Both of the following hold.
\begin{enumerate}
\item Each $E_{2k}$ $(2k\ge4)$ is a rational polynomial of 
$E_{4}, E_{6}$. 
For example, $E_{8}=E_{4}^{2}$, 
$E_{10}=E_{4}E_{6}$, 
$E_{12}=\tff{1}{691}(441E_{4}^{3}+250E_{6}^{2})$, $
E_{14}=E_{4}^{2}E_{6}$.
\item Each $E_{2k}$ $(2k\ge4)$ is a \eh{modular form}  of weight $2k$ on $\text{SL}_{2}(\zz)$ while 
$E_{2}$ is only a quasi-modular form of weight 2 on the same group.
\end{enumerate}

\end{fact}



\begin{thm}[{Ramanujan-Shen's differential equation 
\cite[Theorem D]{shen}}]\label{t2}\mb{}\\
For $m\ge2$, we have 
\begin{eqnarray}
\label{rsdf}
q E'_{2m-2}=
\ff{m-1}{2\p^{2}\z(2m-2)}
\dsum_{k=1}^{m-1}\z(2k)\z(2m-2k)
(E_{2k}E_{2m-2k}-E_{2m}).
\end{eqnarray}
\end{thm}

\begin{rmk}\hf
\begin{enumerate}
	\item Compared to \cite[(23)]{shen}, 
the extra factor $\tf{1}{2}$ shows up in the right hand side since we are taking 
$q=e^{2\piit}$ while $q=e^{\piit}$ therein.
	\item By the symbol $'$, we mean the derivative of a function 
in either a variable $z$ or $q$. 
From context, it will be clear which one $'$ means. In the above theorem, $'=\tff{d}{dq}$. 
\end{enumerate}

\end{rmk}

\begin{ex}
The equations \eqref{rsdf} for $m=2, 3, 4$ are 
\begin{eqnarray}
	qE_{2}'&=\ff{E_{2}^{2}-E_{4}}{12},
	\label{e2d}
	\\qE_{4}'&=\ff{E_{2}E_{4}-E_{6}}{3},
	\label{e4d}
	\\qE_{6}'&=\ff{E_{2}E_{6}-E_{8}}{2}.
	\label{e6d}
\end{eqnarray}
With $E_{8}=E_{4}^{2}$, 
these are nothing but Ramanujan's differential equations.
 For $m=5$, we see that
\begin{eqnarray*}
\label{e8d}
qE_{8}'&=
\ff{2}{3}(E_{2}E_{8}-E_{10})
\end{eqnarray*}
due to the relation %
$E_{4}E_{6}-E_{10}=0$. 
%
\end{ex}

\subsection{Differential equations for Eisenstein series of level 2}
 Our next idea is to consider an analog of \er{ldth1}.
\begin{defn}
For $z\incc, \ta\in\mathbb{H}$, define 
the second kind of Jacobi theta function by 
\[
\th_{2}(z|\ta)=
2
\dsum_{n=0}^{\mug}q^{(2n+1)^{2}/8}\cos(2n+1)z.
\]
\end{defn}

\begin{lem}\label{l1}
Let $f(z)=\th_{2}(z|\ta)$. Then 
\begin{eqnarray}
	\ff{f'}{f}(z)&=
	-\tan z+4
	\dsum_{n=1}^{\mug}\ff{(-1)^{n}q^{n}}{1-q^{n}}\sin 2nz,\label{th2-1}
\end{eqnarray}
\begin{eqnarray}
\ff{f''}{f}(z)&=
	-1+8
	\dsum_{n=1}^{\mug}\ff{nq^{n}}{1-q^{n}}
	+16\dsum_{n=1}^{\mug}\ff{(-1)^{n}q^{n}}{(1-q^{n})^{2}}\cos2nz.\label{th2-2}
\end{eqnarray}
\end{lem}

\begin{proof}
According to \cite[(13), (14)]{shen}, 
if $g(z)=\th_{1}(z|\tau)$, then 
\begin{align*}
	\ff{g'}{g}(z)&=
	\cotz+4
	\dsum_{n=1}^{\mug}\ff{q^{n}}{1-q^{n}}\sin 2nz,
	\\\ff{g''}{g}(z)&=
-1+8
\dsum_{n=1}^{\mug}\ff{nq^{n}}{1-q^{n}}
+16\dsum_{n=1}^{\mug}\ff{q^{n}}{(1-q^{n})^{2}}\cos2nz.
\end{align*}
Note that 
$\th_{2}(z|\tau)=
\th_{1}(z+\tf{\p}{2}|\tau)$ 
and apparently $z\mto z+\tf{\p}{2}$ produces two $(-1)^{n}$s.
\end{proof}

%
\begin{defn}
Let $f(z)=\th_{2}(\piz|\ta)$.
For $k\ge0$, define $G_{2k}^{*}(\ta)$ by 
\begin{eqnarray}\label{g2k*}
\ff{f'}{f}(z)=
-\dsum_{k=0}^{\mug}G_{2k}^{*}(\ta)z^{2k-1}.
\end{eqnarray}
\end{defn}

\begin{prop}
For each $\ta\in\mathbb{H}$, we have 
$G_{0}^{*}(\ta)=0$.
\end{prop}
\begin{proof}

Certainly, $z=0$ is a zero of $\th_{2}'$ since 
\[
\th_{2}'(0|\ta)=-2 \dsum_{n=0}^{\mug}(2n+1)q^{(2n+1)^{2}/8}\sin(2n+1)z
\Bigr|_{z=0}=0
\]
while 
\[
\th_{2}(0|\ta)=
\th_{1}(\tf{\p}{2}|\ta)\ne0.
\]
Thus, $z=0$ is not a pole of $\th_{2}'(\piz|\ta)/\th_{2}(\piz|\ta)$.
\end{proof}

\begin{defn}[signed divisor sum function]
\[
\s^{*}_{s}(n)=-\sum_{d|n}(-1)^{d}d^{s}.
\]
\end{defn}

\begin{lem}
For $k\ge1$, we have 
\[
G_{2k}^{*}(\ta)=
2\cdot 2^{2k}\l(2k)-
2\ff{(2\pii)^{2k}}{(2k-1)!}
\dsum_{N=1}^{\mug}\s^{*}_{2k-1}(N)q^{N}.
\]
\end{lem}
\begin{proof}
%
%
Let $f(z)=\th_{2}(\piz|\ta)$. 
It follows from \eqref{tan}, \eqref{th2-1} that 
\[
\ff{f'}{f}(z)=
\p\ff{\th_{2}'}{\th_{2}}(\piz|\ta)
=
-\p\tan\piz+4\pi
\dsum_{n=1}^{\mug}
\ff{(-1)^{n}q^{n}}{1-q^{n}}\sin 2n\piz
\]
\[=
-2
\dsum_{k=1}^{\mug}
2^{2k}\z^{*}(2k)z^{2k-1}
+
4\pi \dsum_{n=1}^{\mug}
\dsum_{k=1}^{\mug}
\ff{(-1)^{n}q^{n}}{1-q^{n}}
\ff{(-1)^{k-1}}{(2k-1)!}(2n\piz)^{2k-1}
\]
\[=-
\dsum_{k=1}^{\mug}
\k{2\cdot 2^{2k}\l(2k)
-2 \ff{(2\pii)^{2k}}{(2k-1)!}
\k{\dsum_{N=1}^{\mug} \s^{*}_{2k-1}(N)q^{N}
}
}
z^{2k-1}
\]
with interchanging order of the summation, as required.
%
\end{proof}

\begin{defn}[Eisenstein series of level 2]
For $k\ge1$, set 
\[
E_{2k}^{*}(\ta)=\ff{G_{2k}^{*}(\ta)}{2^{2k+1}\l(2k)}
\k{=1-\ff{1}{1-2^{2k}}\ff{4k}{B_{2k}}
\dsum_{n=1}^{\mug}\s_{2k-1}^{*}(n)q^{n}},
\]
$E^{*}_{0}=1$ and 
$\s_{2k-1}^{*}(0)=
\k{-\tff{1}{1-2^{2k}}\tff{4k}{B_{2k}}}^{-1}$.
\end{defn}
\begin{rmk}
$E_{2k}^{*}$ is the same as 
$\mathcal{E}_{2k}$ in {\cite[(1.14)]{hahn2}}
\end{rmk}
\begin{ex}Let us take a look at 
$E_{2k}^{*}$ for $k=1, 2, \ds, 6$.
\begin{align*}
	E_{2}^{*}(\ta)&=
	1+8\dsum_{n=1}^{\mug}\s^{*}_{1}(n)q^{n},
	&E_{4}^{*}(\ta)&=
	1-16\dsum_{n=1}^{\mug}\s^{*}_{3}(n)q^{n},
	\\E_{6}^{*}(\ta)&=
	1+8\dsum_{n=1}^{\mug}\s^{*}_{5}(n)q^{n},
	&E_{8}^{*}(\ta)&=
	1-\ff{32}{17}\dsum_{n=1}^{\mug}\s^{*}_{7}(n)q^{n},
	\\E_{10}^{*}(\ta)&=
	1+\ff{8}{31}\dsum_{n=1}^{\mug}\s^{*}_{9}(n)q^{n},
	&E_{12}^{*}(\ta)&=
	1-\ff{16}{691}\dsum_{n=1}^{\mug}\s^{*}_{11}(n)q^{n}.
\end{align*}
\end{ex}

\begin{lem}\label{l3}
Let $f(z)=\th_{2}(\pi z|\ta)$.
Then
\begin{align*}
	\k{\ff{f'}{f}(z)}^{2}&=
	\k{-\dsum_{k=0}^{\mug}G_{2k}^{*}z^{2k-1}}^{2},
	\\\ff{f''}{f}(z)&=
	-\p^{2}+
8\p^{2}\dsum_{n=1}^{\mug}\ff{nq^{n}}{1-q^{n}}
+16\p^{2}
\dsum_{m=1}^{\mug} \ff{(2\pii)^{2m-2}}{(2m-2)!}
\k{\dsum_{n=1}^{\mug}\ff{(-1)^{n}n^{2m-2}q^{n}}{(1-q^{n})^{2}}}
z^{2m-2},
	\\-
\der{\ff{f'}{f}(z)}&=
\dsum_{m=0}^{\mug}(2m-1)G_{2m}^{*}z^{2m-2}.
\end{align*}
Moreover, $f(z)$ satisfies \eqref{sdf}.
\end{lem}

\begin{proof}
See \er{g2k*} and \er{th2-2}. The last claim follows from that 
\[
f(z)=\th_{2}(\p z|\tau)=
\th_{1}(\p(z+\tf{1}{2})|\tau)
\]
is a solution of Shen's differential equation.
\end{proof}

We are now ready for proving our main theorem.

\begin{thm}\label{t3}
For $m\ge2$, we have 
\begin{eqnarray}\label{e*df}
q{E^{* '}_{2m-2}}=
\ff{2m-2}{\p^{2}\l(2m-2)}
\dsum_{k=1}^{m-1}\l(2k)\l(2m-2k)
(E_{2k}^{*}E_{2m-2k}^{*}-E_{2m}^{*}).
\end{eqnarray}
\end{thm}
\begin{proof}
Lemma \ref{l3} says that 
$f(z)=\th_{2}(\piz|\ta)$ satisfies 
\[
\k{\ff{f'}{f}(z)}^{2}=
\ff{f''}{f}(z)-
\k{\ff{f'}{f}(z)}'.
\]
Equating coefficients of $z^{2m-2}$ $(m\ge2)$ in both sides, we get 
\begin{eqnarray}\label{mid}
\dsum_{k=0}^{m}G_{2k}^{*}
G_{2m-2k}^{*}
=
16\p^{2}
 \ff{(2\pii)^{2m-2}}{(2m-2)!}
\dsum_{n=1}^{\mug}\ff{(-1)^{n}n^{2m-2}q^{n}}{(1-q^{n})^{2}}
+(2m-1)G_{2m}^{*}.
\end{eqnarray}
By differentiating 
\[
\dsum_{n=1}^{\mug}\ff{(-1)^{n}n^{2m-3}q^{n}}{1-q^{n}}
=
\dsum_{N=1}^{\mug}\sum_{n|N}(-1)^{n}n^{2m-3}q^{N}
=-
\dsum_{N=1}^{\mug}\s_{2m-3}^{*}(N)q^{N}
\]
and multiply by $q$, we get another expression of 
the sum in the middle of \er{mid}:
\[
\dsum_{n=1}^{\mug}\ff{(-1)^{n}n^{2m-2}q^{n}}{(1-q^{n})^{2}}=
-
\dsum_{N=1}^{\mug}N\s_{2m-3}^{*}(N)q^{N}.
\]
Hence 
\begin{align*}
	q 
E_{2m-2}^{*'}&=
\ff{1}{2^{2m-1}\l(2m-2)}\,
\ff{-2(2\pii)^{2m-2}}{(2m-3)!}
\dsum_{N=1}^{\mug}N\s^{*}_{2m-3}(N)q^{N}
	\\&=\ff{2(2\pii)^{2m-2}}{2^{2m-1}\l(2m-2)(2m-3)!}
\dsum_{n=1}^{\mug}
\ff{(-1)^{n}n^{2m-2}q^{n}}{(1-q^{n})^{2}}
	\\&=\ff{2(2\pii)^{2m-2}}{2^{2m-1}\l(2m-2)(2m-3)!}
	\\&\ti
\ff{(2m-2)!}{16\p^{2}(2\pii)^{2m-2}}
\k{\dsum_{k=0}^{m}G_{2k}^{*}
G_{2m-2k}^{*}-
{(2m-1)G_{2m}^{*}}
}
	\\&=\ff{2m-2}{\l(2m-2)\p^{2}}
{\dsum_{k=1}^{m-1}
\l(2k)\l(2m-2k)
(E_{2k}^{*}
E_{2m-2k}^{*}-E_{2m}^{*})
}
\end{align*}
as we used \eqref{z*2m} for the last equality.
\end{proof}

\begin{ex}
For $m=2, 3, 4$, we have 
\begin{align}
	q E_{2}^{*'}&
=\ff{1}{4}(E_{2}^{* 2}-E_{4}^{*}),
\label{e2*d}
\\q E_{4}^{*'}&=
E_{2}^{*}E_{4}^{*}-E_{6}^{*},\label{e4*d}
\\q E_{6}^{*'}&=
\ff{1}{8}
\k{12E_{2}^{*}E_{6}^{*}
+5E_{4}^{* 2}-{17}E_{8}^{*}}.
\label{e6*d}
\end{align}
We can continue this to get
\[
q E_{8}^{*'}=
\ff{1}{17}
\k{34E_{2}^{*}E_{8}^{*}
+28E_{4}^{*}E_{6}^{*}-{62}E_{10}^{*}}.
\]
\end{ex}
For notational simplicity, let us write 
\[
A=E_{2}^{*}, \q B=E_{4}^{*}, \q
C=\ff{E_{6}^{*}}{E_{4}^{*}}.
\]
\begin{prop}[{\cite[(2.9)]{hahn2}}]
We have $E_{8}^{*}=\tff{1}{17}(9B^{2}+8BC^{2}).$
\end{prop}

As a consequence of this and \eqref{e6*d}, 
we find the polynomial expression 
\begin{eqnarray}\label{e6*d2}
qE_{6}^{*'}=
\ff{1}{2}\k{3ABC-B^{2}-2BC^{2}}
\end{eqnarray}
in terms of only $A, B, C$.
\begin{rmk} 
Hahn \cite{hahn} (and other authors) proved that the 
triple $(\P, \E, \Q)=(A, C, B)$ forms the system
\[
q\P'=\ff{\P^{2}-\Q}{4},\q 
q\E'=\ff{\P\E-\Q}{2},\q 
q\Q'=\P\Q-\E\Q.
\]
We can also derive this from Theorem \ref{t3}. 
The first and third equations are \eqref{e2*d} and \eqref{e4*d}. 
For the second one, observe that 
\begin{align*}
	q\E'&=q
\k{\ff{E_{6}^{*}}{E_{4}^{*}}}'
	\\&=\ff{\tf{1}{2}\k{3ABC-B^{2}-2BC^{2}}B-BC(AB-BC)
}{B^{2}}
	\\&=\ff{1}{2}(AC-B)=\ff{\P\E-\Q}{2}.
\end{align*}
\end{rmk}

\section{Arithmetic applications}

\renewcommand{\d}{\Delta}

\subsection{Divisor sum functions}

Ramanujan's differential equations imply 
certain arithmetic identities, for example, 
\[
\s_{3}(n)=
\ff{6}{5}
\k{n\s(n)+
2\dsum_{j=0}^{n}\s(j)\s(n-j)}
\]
for all $n\ge0$. Ramanujan-Shen differential equation  for $m=7$ is particularly simple due to the relations 
$E_{4}E_{10}-E_{14}=0$, $E_{6}E_{8}-E_{14}=0$ so that 
\[
qE'_{12}=
E_{2}E_{12}-E_{14}.\]
This differential equation with 
\[
E_{12}=
1+\ff{65520}{691}
\dsum_{n=1}^{\mug}\s_{11}(n)q^{n}, \q
E_{14}=
1-24
\dsum_{n=1}^{\mug}\s_{13}(n)q^{n}, 
\]
shows the following formula.
\begin{cor}
Recall that $\s(0)=-\tf{1}{24}, 
\s_{11}(0)=\tf{691}{65520}$ and 
$\s_{13}(0)=-\tf{1}{24}$. 
\label{t7}
For all $n\ge0$, we have 
\[
\s_{13}(n)=
\ff{2730}{691}
\k{24
\sum_{j=0}^{n}\s(j)\s_{11}(n-j)+n\s_{11}(n)
}.
\]
\end{cor}
See Hahn \cite{hahn}, Huard-Ou-Spearman-Williams \cite{hosw}, for other examples. Let us prove an analogous identity from our differential equations. 
Recall that $\s^{*}(0)=\tf{1}{8}$ and 
$\s^{*}_{3}(0)=-\tf{1}{16}$.
\begin{thm}\label{t5}
For $n\ge0$, we have 
\begin{align*}
	\s_{3}^{*}(n)
&=2n\s^{*}(n)-4
\dsum_{j=0}^{n}\s^{*}(j)\s^{*}(n-j).
\end{align*}
\end{thm}
\begin{proof}
Equate the coefficients of $q E_{2}^{*'}=
\tf{1}{4}(E_{2}^{* 2}-E_{4}^{*}).$
%
\end{proof}
 
In this way, we can derive infinitely many more such formulas from \er{e*df}. 

\subsection{Ramanujan's tau function}

Recall that the \eh{Ramanujan's discriminant} is 
\[
\d(q)=
q\dprod_{n=1}^{\mug} (1-q^{n})^{24}.
\]
Let $\tau(n)$ be the \eh{Ramanujan's tau function}, i.e.,  \[
\d=
\dsum_{n=1}^{\mug}\tau(n)q^{n}.
\]
In fact, $\ta(1)=1, \ta(2)=-24, \ta(3)=252, \ta(4)=-1472$. Set $\ta(0)=0$. 
Many authors proved congruence properties of $\ta(n)$ such as 
\[
\hou{\ta(n)}{\s_{11}(n)}{691}.
\]
We will show new formulas in Corollaries \ref{c1} and \ref{c2}. 
For the moment, let us review some history because of its importance in number theory.

\begin{thm}[Mordell 1917 \cite{mordell}]
Both of the following hold.
\begin{enumerate}
	\item $\ta$ is multiplicative. 
	\item $\ta(p^{k+1})=\tau(p)\ta(p^{k})-p^{11}\ta(p^{k-1})$ for all prime $p$ and $k\ge1$.
\end{enumerate}
\end{thm}

\begin{cj}[Ramanujan 1916 \cite{ramanujan}] 
For all prime $p$, $|\tau(p)|\le 2p^{11/2}$.
\end{cj}
Deligne \cite{de} proved this in 1974.

\begin{cj}[Lehmer 1947 \cite{le}]
$\ta(n)\ne$ 0 for $n\ge1$.
\end{cj}
This conjecture remains open. However, researchers have been making progress on its variants.

\begin{thm}[M.R.Murty-V.K.Murty-Shorey 1987 \cite{mms}]
For each odd integer $a$, there exists only finitely many $n\innn$ such that 
$\ta(n)=a$.
\end{thm}
For more recent progress on ``forbidden" $\tau$-values, refer to Balakrishnan-Ono-Tsai \cite{bot}, 
Bennett-Gherga-Patel-Siksek \cite{bgps}, 
Lakein-Larsen \cite{ll} and Moree \cite{mor}.

In the theory of $j$-invariant, we often encounter the following expression.
 
\begin{fact}
The following relation holds.
\begin{eqnarray}\label{dis}
\d=\ff{{E_{4}}^{3}-{E_{6}}^{2}}{1728}.
\end{eqnarray}
\end{fact}
Note that we can express it with a determinant as 
\[
\d=
\ff{1}{1728}
\left|\begin{array}{cc}
	E_{4}&E_{6}   \\
	E_{6}&E_{8}   \\
\end{array}\right|.
\]
More generally, 
Milne studied Hankel determinants of $(E_{2k})_{2k\ge4}$.
\begin{fact}[Garvan's identity {\cite[Theorem 1.4]{milne}}]We have 
\[
\left|\begin{array}{ccc}
	E_{4}&   E_{6}&E_{8}   \\
	E_{6}&   E_{8}&E_{10}   \\
	E_{8}&   E_{10}&E_{12}
	\end{array}\right|
=-
\ff{250}{691} (1728\d)^{2}.
\]
\end{fact}
Milne \cite{milne} computed other such determinants which all involve $\d$.

Here let us try computing minors of 
$|E_{2(i+j)-4}|_{i, j\ge1}^{\mg}$ including $E_{0}=1$ and $E_{2}$.
Certainly, some of them are a multiple of $qE_{2k}'$  such as 
\[
\left|\begin{array}{cc}
	E_{0}&E_{2}   \\
	E_{2}&E_{4}   \\
\end{array}\right|=-12qE_{2}', \q 
\left|\begin{array}{cc}
	E_{0}&E_{2}   \\
	E_{4}&E_{6}   \\
\end{array}\right|=-3qE_{4}',\]
\[
\left|\begin{array}{cc}
	E_{2}&E_{4}   \\
	E_{4}&E_{6}   \\
\end{array}\right|=2qE_{6}', \q 
\left|\begin{array}{cc}
	E_{2}&E_{6}   \\
	E_{4}&E_{8}   \\
\end{array}\right|=\ff{3}{2}qE_{8}'.
\]
Now it is natural to reinterpret \eqref{dis} as a  differential equation.
%
%


\begin{lem}\label{l4}
We have \[
1728\d=3E_{6}qE_{4}'-2E_{4}qE_{6}'.
\]
\end{lem}
\begin{proof}
\[
\left|\begin{array}{cc}
	E_{4}&E_{6}   \\
	E_{6}&E_{8}   \\
\end{array}\right|=
\left|\begin{array}{cc}
	E_{4}&E_{6}   \\
	E_{6}-E_{2}E_{4}&E_{8}-E_{2}E_{6}
	\\
\end{array}\right|=
3E_{6}qE_{4}'-2E_{4}qE_{6}'.
\]
\end{proof}

Recall that 
$\s_{3}(0)=\tf{1}{240}$ and 
$\s_{5}(0)=-\tf{1}{504}$.
\begin{thm}\label{t8}
For $n\ge0$, 
\begin{align*}
	\ta(n)&=70
	\sum_{
	\substack{j, k\ge0\\j+k=n	
	}
	}(2k-3j)\s_{3}(j)\s_{5}(k).
\end{align*}
\end{thm}

\begin{proof}
Use Lemma \ref{l4}.
\end{proof}

\begin{cor}\label{c1}
For all $n\ge0$, 
\[
\hou{\ta(n)}{\ff{n}{12}(5\s_{3}(n)+7\s_{5}(n))}{70}.
\]
\end{cor}
\begin{proof}We have 
\[
\ta(n)=
70
\k{\ff{2n}{240}\s_{5}(n)-3n\ff{\s_{3}(n)}{-504}}+70
\sum_{
	\substack{j, k>0\\j+k=n	
	}
	}(2k-3j)\s_{3}(j)\s_{5}(k)
\]
and the sum in the right hand side is an integer.
\end{proof}
\begin{ex}Observe that 
\begin{align*}
	\ta(3)-\ff{3}{12}(5\s_{3}(3)+7\s_{5}(3))&=
	252-\ff{1}{4}(5\cdot 28+7\cdot 244)
=\hou{-210}{0}{70}.
\end{align*}
\end{ex}

\begin{thm}\label{t314}
We have 
\[
\ta(n)=
2
\sum_{
	\substack{j, k\ge0\\j+k=n	
	}
	}(3j-2k)\s_{3}^{*}(j)\s_{5}^{*}(k).
\]
\end{thm}
\begin{proof}
Using $E_{2}=3A-2C$ \cite[Lemma 2.5]{hahn2} 
and \er{e2d}, \er{e4d}, \er{e6d}, \er{e2*d}, \er{e4*d}, \er{e6*d2}, 
we have 
\[
E_{4}=E_{2}^{2}-12qE_{2}'=-3B+4C^{2}
\]
and moreover 
\[
E_{6}=E_{2}E_{4}-3qE_{4}'=
9BC-8C^{3}
\]
with some algebra. It turns out that 
\begin{eqnarray}\label{delta}
\d=\ff{E_{4}^{3}-E_{6}^{2}}{1728}=
\ff{
(-3B+4C^{2})^{3}-(9BC-8C^{3})^{2}
}{1728}
=
-\ff{{E_{4}^{*}}^{3}-{E_{6}^{*}}^{2}}{64}.
\end{eqnarray}
Therefore, with \er{e4*d}, \er{e6*d2}, 
we will arrive at
\[
-64\d=3E_{6}^{*}qE_{4}^{*'}-2E_{4}^{*}qE_{6}^{*'}.
\]
\end{proof}
\begin{cor}\label{c2}We have 
\[
\hou{\ta(n)}{\ff{n}{4}(3\s_{3}^{*}(n)+\s_{5}^{*}(n))}{2}.
\]
As a consequence, 
\[
\mb{$\ta(n)$ is odd} \iff 
\hou{n(3\s_{3}^{*}(n)+\s_{5}^{*}(n))}{4}{8}.
\]
\end{cor}

\section{Modular forms}

\subsection{Serre Derivative}

Next, we develop the discussions on {modular forms} (without proofs). 
Let 
\[
\Gam=\te{SL}_{2}(\zz)=\Set{
\begin{pmatrix}
	a&b   \\
	c&d   
\end{pmatrix}\in M_{2}(\zz)}{ad-bc=1}
\]
be the full modular group
and 
\[
\Gam_{0}(2)=\Set{\begin{pmatrix}
	a&b   \\
	c&d   
\end{pmatrix}\in \te{SL}_{2}(\zz)}{
\hou{c}{0}{2}}
\]
its congruence subgroup of level 2.
Denote by $M_{2k}=M_{2k}(\Gam_{0}(2))$ 
the set of modular forms of weight $2k$. This is a finite dimensional $\cc$-vector space. 
In particular, 
$\dim M_{2}=1$ and 
\[
\ff{E_{6}^{*}}{E_{4}^{*}}=
1+24\sum_{n=1}^{\mg}
\sh(n)q^{n}, \q 
\sh(n)=\sum_{\substack{d|n\\d \text{ odd}}}d
\]
forms its basis. See \mbox{\cite[p. 10, Proposition 3.4. 1 and p.13 (3.8a)]{maier}}; our $\tff{E_{6}^{*}}{E_{4}^{*}}$ coincides with Maier's $\mathcal{A}_{4}^{2}$.
 Note that $E_{2}^{*}$ is \eh{not} a 
modular, but a quasi-modular form of weight 2 on 
$\Gam_{0}(2)$. 
Recall our simplified notation
\[
A=E_{2}^{*}, \q B=E_{4}^{*}, \q
C=\ff{E_{6}^{*}}{E_{4}^{*}}.
\]
The weight of a monomial 
$A^{a}B^{b}C^{c}$
is $2a+4b+2c$. 
By $\cc[A, B, C]_{2k}$ we mean the set of 
all complex homogenous polynomials in $A, B, C$ of weight $2k$ and $0$. We refer any $f\in \cc[A, B, C]$ 
as a quasi-modular form on $\Gam_{0}(2)$.
\begin{prop}[{\cite[Proposition 3.5]{maier}}]\label{m2k}
If $k\ge1$, then 
\[
M_{2k}(\Gam_{0}(2))=
\cc
\left[B, C\right]_{2k}.
\]
Moreover, $\dim_{\cc} M_{2k}=[\tf{2k}{4}]+1$ 
and $\{B^{j}C^{k-2j}\mid j\in\zz, 0\le j\le \tf{k}{2}\}$ forms its basis.
\end{prop}

{\renewcommand{\arraystretch}{1.25}
\begin{table}[h!]
\caption{vector spaces $M_{2k}, S_{2k}$}
\label{tb1}
\begin{center}
\begin{tabular}{@{} c|c|c||c|c|c @{}}
$M_{2k}$&$\dim M_{2k}$&basis&
$S_{2k}$&$\dim S_{2k}$ &basis \\\hline
$M_{2}$&1&$C$&$S_{2}$&0&---  \\\hline
$M_{4}$&2&$B, C^{2}$&$S_{4}$&0&---    \\\hline
$M_{6}$&2&$BC, C^{3}$&$S_{6}$&0&---    \\\hline
$M_{8}$&3&$B^{2}, BC^{2}, C^{4}$&
$S_{8}$&1&$BD$  \\\hline
$M_{10}$&3&$B^{2}C, BC^{3}, C^{5}$&
$S_{10}$&1&$BCD$  \\\hline
$M_{12}$&4&$B^{3}, B^{2}C^{2}, BC^{4}, C^{6}$&
$S_{12}$&2&$\d=B^{2}D, BD^{2}$  \\\hline
\end{tabular}
\end{center}
\end{table}}

%

\begin{ex}[{\cite[(2.10), (2.11)]{hahn2}}]
\label{p1}Observe that 
\begin{align*}
	E_{10}^{*}&=\ff{1}{31}(27B^{2}C+4BC^{3}) \,\,(\in M_{10}),
	\\E_{12}^{*}&=\ff{1}{691}(
	189B^{3}+486B^{2}C^{2}+16BC^{4}) \,\,(\in M_{12}).\label{e12*d}
\end{align*}
\end{ex}
We expect that such coefficients for $E_{2k}^{*}$ 
be all positive  rational numbers. In fact, this is true as shown in Theorem \ref{t49}. 
For this purpose, we introduce the following operator.
%



\begin{defn}
For $f\in \cc
\left[A, B, C\right]_{2k}$, its \eh{Serre derivative} is 
\[
\del_{2k}f=qf'-\ff{2k}{4}Af.
\]
\end{defn}
Whenever we understand the weight of $f$, we just write $\del f$.

\begin{rmk}\hf
\begin{enumerate}
	\item The coefficient of 
	$E_{2}^{*}E_{2m-2}^{*}$ 
	in right hand side of 
	\er{e*df} happens to be 
\[\ff{2m-2}{\p^{2}\z^{*}(2m-2)}
2\z^{*}(2)\z^{*}(2m-2)=\ff{2m-2}{4}.
\]
So both sides of \er{e*df} includes $\del E_{2m-2}^{*}$ implicitly.
\item Similarly, if we define the Serre derivative 
of $f\in \cc[E_{2},E_{4}, E_{6}]_{2m-2}$ by 
\[
\partial f=qf'-\ff{2m-2}{12}E_{2}f,
\]
the coefficient of $E_{2}E_{2m-2}$ in right hand side of \er{rsdf} is $\tff{2m-2}{12}$.
\end{enumerate}
\end{rmk}

Now we can show the new closed system of differential equations as mentioned in Introduction.
\begin{prop}\label{p4}
The following system of differential equation holds.
\[
\del A=-\ff{1}{4}(A^{2}+B), \q 
\del B=-BC, \q 
\del C=-\ff{1}{2}B.
\]
As a consequence, if $k\ge1$ and 
$f \in \cc[A, B, C]_{2k}$, then 
$\del f \in \cc[A, B, C]_{2k+2}$ .
\end{prop}
\begin{proof}
\begin{align*}
	\del A&=qA'-\ff{2}{4}A^{2}=-\ff{1}{4}(A^{2}+B),
	\\\del B&=qB'-\ff{4}{4}AB=-BC,
	\\\del C&=q\der{\ff{E_{6}^{*}}{E_{4}^{*}}}-\ff{2}{4}AC
	\\&=
\ff{\tf{1}{2}\k{3ABC-B^{2}-2BC^{2}}B-BC(AB-BC)
}{B^{2}}
-\ff{1}{2}AC=
-\ff{1}{2}B.
\end{align*}
\end{proof}

\begin{cor}
If $k\ge1$, $f\in M_{2k}$, then $\del f\in M_{2k+2}$.
\end{cor}
The operator $\del$ satisfies Leibniz rule in the sense that 
for $f\in M_{2k}, g\in M_{2l},$ we have 
\[
\del_{2k+2l}(fg)=\del_{2k}(f)g+f\del_{2l}(g).
\]
With Proposition \ref{p4}, we can inductively compute $\del f$ for any $f\in \cc[A, B, C]$.
\begin{ex}Observe that 
\begin{align*}
	\del C^{2}&=2C\del(C)=-BC.
	\\\del (BC)&=\del(B)C+B\del(C)
=-\ff{1}{2}B^{2}-BC^{2}.
	\\\del C^{3}&=
(\del C)C^{2}+C(\del C^{2})=
-\ff{3}{2}BC^{2}.
\end{align*}
\end{ex}

Let 
$\qq_{+}$ $(\qq_{-})$ denote the set of all positive (negative) rational numbers.
Further, denote $\qq_{+}[x_{1}, \ds, x_{n}]$
 (resp. $\qq_{-}[x_{1}, \ds, x_{n}]$) 
the set of polynomials in $x_{1}, \ds, x_{n}$ with all coefficients positive (resp. negative) rational numbers.

\begin{lem}\label{neg}Let $k\ge1$. 
For an integer $0\le j\le \tff{k}{2}$, 
$\del(B^{j}C^{k-2j})\in B\qq_{-}\left[B, C\right]$.
\end{lem}

\begin{proof}
First, by induction on $n$, we can show that 
for each $n\in\nn$, 
\[
\del B^{n}, \del C^{n}\in B\qq_{-}[B, C].
\]
Now let $j, k$ be integers such that $k\ge1$ and 
$0\le j\le \tff{k}{2}$. Since 
$\del(B^{j}C^{k-2j})=
\del(B^{j})C^{k-2j}+B^{j}\del(C^{k-2j})$, 
we have 
$\del(B^{j}C^{k-2j})\in B\qq_{-}\left[B, C\right]$.
\end{proof}

\begin{thm}\label{t49}
For $m\ge2$, we have 
\[
E_{2m}^{*}\in B\qq_{+}
\left[B, C\right].
\]
\end{thm}
\begin{proof}
Induction on $m$.
If $m=2$, then 
$E_{4}^{*}=B\in B\qq_{+}[B, C]$ itself.
Suppose $m\ge3$. 
Now, \er{e*df} together with \er{z*2m} imply
\[
E_{2m}^{*}=
\ff{1}{\al_{2m}}
\k{
\ff{2m-2}{\p^{2}\l(2m-2)}
\dsum_{k=2}^{m-2}\l(2k)\l(2m-2k)
E_{2k}^{*}E_{2m-2k}^{*}
-\del E_{2m-2}^{*}
}
\]
with 
\[
\al_{2m}=\ff{2m-2}{\p^{2}\l(2m-2)}\ff{2m-1}{2}\z^{*}(2m)\in \qq_{+}.
\]
We are assuming for our induction that 
\[
E_{4}^{*}, \ds, 
E_{2m-4}^{*}, E_{2m-2}^{*}\in B\qq_{+}[B, C]
\]
and hence 
by Proposition \ref{m2k} and Lemma \ref{neg}, 
we have $-\del E_{2m-2}^{*}
\in B\qq_{+}[B, C]$.
Conclude that 
\[
E_{2m}^{*}\in B\qq_{+}
\left[B, C\right].
\]
%
%
\end{proof}
\begin{cor}
For all $m\ge2$, we have $E_{2m}^{*}\in 
B\qq_{+}(B, C)_{2m-4}$.
\end{cor}
\begin{proof}
Thanks to Theorem \ref{t49}, 
$E^{*}_{2m}=Bf(B, C)$ 
for some $f\in \mathbb{Q}_{+}[B, C]$ (and indeed $f\ne0$). Further, Proposition \ref{m2k} says that 
$E_{2m}^{*}$ is a polynomial in $B, C$ of weight $2m$. 
Thus, $f$ must have weight $2m-4$ and hence 
\[
E^{*}_{2m}=Bf(B, C)\in B\qq_{+}(B, C)_{2m-4}.
\]
\end{proof}
%

\subsection{$D$ and Ramanujan's discriminant $\d$}

The map $\del_{2k}:M_{2k}\to M_{2k+2}$ is linear. 
For understanding its kernel, it is useful to consider the following form.
\begin{defn}
$D=-\tf{B-C^{2}}{64}$.
\end{defn}
Observe that $D\in M_{4}=\cc B\oplus \cc C^{2}
$ and it contains monomials with \eh{both} positive and negative coefficients.
\begin{prop}We have $\del D=0$.
\end{prop}
Thus, $\te{ker}\, \del_{4}=\cc D$ is a 1-dimensional subspace of $M_{4}$.
With all results we have shown, it is easy to check that the following family of differential equations holds for $\del=\del_{4}$:
\begin{align*}
\del \ff{E_{6}^{*2}}{E_{4}^{*2}}=
	\del E_{4}^{*}&=\del \ff{E^{*}_{8}}{E^{*}_{4}}=
	\del \ff{E^{*}_{10}}{E^{*}_{6}}=
	\del E_{4}.
\end{align*}
Thus, for example, $E_{4}^{*}-\tff{E^{*}_{8}}{E^{*}_{4}}\,(\in M_{4})$ must be a multiple of $D$ which immediately implies the following.
\begin{lem}\label{l5}
We have 
\[
\left|\begin{array}{cc}
	E_{0}^{*}&E_{4}^{*}   \\
	E_{4}^{*}&E_{8}^{*}   \\
\end{array}\right|=
\ff{512}{17}BD.
\]
\end{lem}

\begin{rmk}\hf
\begin{enumerate}
\item 
Let $S_{2k}$ denote the set of \eh{cusp forms} of weight $2k$ on 
$\Gam_{0}(2)$. We have 
$\dim S_{2}=0$ and for $k\ge2$, $\dim S_{2k}=t-1$ where $t=[\tf{2k}{4}]$ (integer part).
In particular, $\dim S_{8}=1$ and $S_{8}=\cc BD$ as in Table \ref{tb1}. 
If 
$\hou{2k}{0}{4}$, 
\[
BD^{t-1}, B^{2}D^{t-2}, \ds, B^{t-1}D
\]
forms a basis of $S_{2k}$.
If 
$\hou{2k}{2}{4}$, 
\[
BCD^{t-1}, B^{2}CD^{t-2}, \ds, B^{t-1}CD
\]
forms a basis of $S_{2k}$ as Imamo$\bar{\text{g}}$lu-Kohnen discussed \cite[p.828]{ik}. We remark that 
$B=E_{4}^{i\mug}, D=E_{4}^{0}$ and  
$BC=E_{6}^{i\mug}$ in their notation.
\item Quite similarly, $\partial \d=0$ 
and thus $\partial f=0 \iff f\in \cc \d$ for 
$f\in M_{12}(\text{SL}_{2}(\zz))$.
Brundaban-Ramakrishnan \cite{br} 
found a series of formulas for $\d$ based on this idea.
In the literature, authors often discuss these equations  equivalently as $q\d'=E_{2}\d$ and $qD'=E_{2}^{*}D$. 
	\item Coefficients of $D$ has also a nice arithmetical interpretation \cite[Theorem 2.4 (2.51)]{hahn}:  
\[
D=
q 
\dsum_{n=0}^{\mug}\delta_{8}(n)q^{n}
\]
where $\delta_{8}(n)$ is the number of representations of $n$ as the sum of 8 triangular numbers. 
The problem of finding formulas for 
$\delta_{s}(n)$ is one of the classical topics in number theory.
\end{enumerate}
\end{rmk}

\subsection{Sum of squares}

Let 
$\th_{3}(q)=\sum_{m\in \zz}q^{m^{2}}$ be 
the third kind of Jacobi theta function. 
For $s\ge1$, let $r_{s}(n)\in \zz_{\ge0}$ such that 
\[
\th_{3}^{s}(q)=
\dsum_{n=0}^{\mug}r_{s}(n)q^{n},
\]
that is, 
\[
r_{s}(n)=|\{(m_{1}, \ds, m_{s})\in \zz^{s}\mid 
m_{1}^{2}+\cd+m_{s}^{2}=n\}|.
\]
There is long history on these numbers; Euler, Fermat, Gauss, Hardy, Legendre, Mordell, Rankin etc. discovered  formulas for $r_{2k}$. In  particular, Lagrange proved $r_{4}(n)>0.$ 
Jacobi further proved that 
\begin{align*}
	r_{2}(n)&=4
\k{
\sum_{
\substack{d|n\\
d\equiv 1 \text{ mod }4}
}d-
\sum_{
\substack{d|n\\
d\equiv 3 \text{ mod }4
}
}d
}
, 
	\\r_{4}(n)&=8
\sum_{
\substack{d|n\\
d\not\equiv 0 \text{ mod }4
}
}d,
	\\r_{6}(n)&=4
\sum_{
\substack{d|n\\
d\text{ odd}}
}(-1)^{(d-1)/2}
\k{\k{\ff{2n}{d}}^{2}-d^{2}
},
	\\r_{8}(n)&=16(-1)^{n}
\sum_{
\substack{d|n}
}(-1)^{d}d^{3}.
\end{align*}
In particular, the last identity 
implies 
\begin{eqnarray}\label{theta}
E_{4}^{*}(q)=\th_{3}^{8}(-q).
\end{eqnarray}
More recently, researchers discovered 
formulas of $r_{2k}(n)$ for $2k=16, 24, 32$ as in  \cite{chanchua, ik, milne2}. We add some more to the list.
\begin{thm}We have\label{t9}
\[
r_{16}(n)=
(-1)^{n}\ff{32}{17}
\k{256
\dsum_{j=0}^{n-1}\s^{*}_{3}(j)\delta_{8}(n-j-1)-\s_{7}^{*}(n)
}.
\]
\end{thm}
\begin{proof}
Lemma \ref{l5} with $E_{4}^{*}(q)=\th_{3}^{8}(-q)$
shows that 
$\th_{3}^{16}(-q)=E_{4}^{*2}(q)=E_{8}^{*}(q)-\tff{512}{17}BD$. 
Equating coefficients of $q^{n}$ $(n\ge1)$ in left and right hand sides, we have 
\[
(-1)^{n}r_{16}(n)=
-\ff{32}{17}\s^{*}_{7}(n)-\ff{512}{17}
\sum_{
\substack
{j, k\ge0\\
j+k=n
}
}
(-16\s_{3}^{*}(j))
\delta_{8}(k-1).
\]
This is equivalent to our assertion.
\end{proof}

\begin{fact}[Ramanujan's 24-square theorem]
\[
r_{24}(n)=
\ff{1}{691}
\left(
16\s_{11}(n)-32\s_{11}(\tf{n}{2})
+65536\s_{11}(\tf{n}{4})
\right.
\]
\[
\left.
+33152(-1)^{n-1}\tau(n)-
65536\tau(\tf{n}{2})
\right)
\]
where $\s_{11}(x)=\ta(x)=0$ if $x\not\in \zz_{\ge0}$.
\end{fact}
Obviously, this comes from the calculus of $\th_{3},  E_{12}$ and $\d$. Let us figure out some relation among $\th_{3}, E_{2k}^{*}$ and $\d$. We then come across certain coincidence of coefficients so that we find  a nice relation between $r_{24}(n)$ and $\ta(n)$ with $(-1)^{n}$.
\begin{thm}\label{t10}We have 
\begin{align*}
	r_{24}(n)&=(-1)^{n}64
\k{
\dsum_{j=0}^{n}\s_{5}^{*}(j)\s_{5}^{*}(n-j)-\tau(n)
}
	\\&=
(-1)^{n}\ff{512}{17}
\k{
\dsum_{j=0}^{n}\s_{3}^{*}(j)\s_{7}^{*}(n-j)-\tau(n)
}.
\end{align*}
\end{thm}
\begin{proof}
The first equality follows from 
$\th_{3}^{24}(-q)=E_{4}^{*3}(q)=E_{6}^{*2}(q)-64\d(q)$ (See \eqref{theta}, \eqref{delta}), that is, 
\[
(-1)^{n}r_{24}(n)=
\dsum_{j=0}^{n}8^{2}\s_{5}^{*}(j)\s_{5}^{*}(n-j)-64\ta(n)
\]
for all $n\ge0$. 
For the second equality, use Lemma \ref{l5}.
\end{proof}
\begin{cor}\label{c10}
For $n\ge0$, the following are equivalent:
\begin{quote}
	\begin{enumerate}
		\item $n$ is odd.
		\item $\tau(n)> \dsum_{j=0}^{n}\s_{5}^{*}(j)\s_{5}^{*}(n-j).$
		\item $\tau(n)> \dsum_{j=0}^{n}\s_{3}^{*}(j)\s_{7}^{*}(n-j).$
	\item 
	$\dsum_{j=0}^{n}\s_{5}^{*}(j)\s_{5}^{*}(n-j)>\dsum_{j=0}^{n}\s_{3}^{*}(j)\s_{7}^{*}(n-j)$.
	\end{enumerate}
\end{quote}
\end{cor}
\begin{proof}
First, we know $r_{24}(n)\ge r_{4}(n)>0$.
Second, Theorem \ref{t10} further says
\[
9r_{24}(n)=(-1)^{n}512
\k{\dsum_{j=0}^{n}\s_{5}^{*}(j)\s_{5}^{*}(n-j)
-
\dsum_{j=0}^{n}\s_{3}^{*}(j)\s_{7}^{*}(n-j)
}.
\]
Now we see that (i), (ii), (iii), (iv) are all equivalent.
\end{proof}
{\renewcommand{\arraystretch}{1.25}
\begin{table}[h!]
\label{tb2}
\caption{values of arithmetic functions}
\begin{center}
\begin{tabular}{c|ccccccccccccccccccc}
$n$&0&1&2&3&4\\\hline
$\s_{3}^{*}(n)$&$-\tf{1}{16}$&1&$-7$&28&$-71$\\\h
$\s_{5}^{*}(n)$&$\tf{1}{8}$&1&$-31$&244&$-1055$\\\h
$\s_{7}^{*}(n)$&$-\tf{17}{32}$&1&$-127$&2188&$-16511$\\\h
$
\sum_{j=0}^{n}\s_{3}^{*}(j)\s_{7}^{*}(n-j)
$
&$\tf{12}{517}$&$-\tf{19}{32}$&$\tf{405}{32}$&
$-\tf{2285}{8}$&$\tf{133589}{32}$\\\h
$
\sum_{j=0}^{n}\s_{5}^{*}(j)\s_{5}^{*}(n-j)
$
&$\tf{1}{64}$&$\tf{1}{4}$&$\tf{33}{32}$&
$-1$
&$\tf{37928}{32}$\\\h
$\ta(n)$&0&1&$-24$&252&$-1472$\\\h
\end{tabular}
\end{center}
\end{table}%
}

%


It may be of interest to investigate a system of 
minors in $(E_{2(i+j)}^{*})$ as Garvan and Milne studied. 
Here, we finish recording several examples.
\begin{align*}
	\left|\begin{array}{cc}
	E_{4}^{*}&E_{6}^{*}   \\
	E_{6}^{*}&E_{8}^{*}   \\
\end{array}\right|
&=-\ff{2^{6}3^{2}}{17}\d,
\\
	\left|\begin{array}{cc}
	E_{4}^{*}&E_{8}^{*}   \\
	E_{6}^{*}&E_{10}^{*}   \\
\end{array}\right|&=
-\ff{2^{8}3^{2}5}{17\cdot 31}C\d,
	\\\left|\begin{array}{cc}
	E_{6}^{*}&E_{8}^{*}   \\
	E_{8}^{*}&E_{10}^{*}   \\
\end{array}\right|
&=\ff{2^{6}3^{2}}{17^{2}31}
(279B-92C^2)\d
\end{align*}
and 
\[
\left|\begin{array}{ccc}
	E_{4}^{*}&   E_{6}^{*}&E_{8}^{*}   \\
	E_{6}^{*}&   E_{8}^{*}&E_{10}^{*}   \\
	E_{8}^{*}&   E_{10}^{*}&E_{12}^{*}
	\end{array}\right|
=
-\ff{2^{13}3^{5}5^{2}}{17^{3}31^{2}691}
(961B+3136C^{2})BD\d.
\]
%


\paragraph*{\bf Acknowledgment.}
This work is based on the author's talk in Hiroshima-Sendai Number Theory Workshop at Hiroshima  University, Japan, in July 2023. He thanks its organizers and the audience who gave valuable comments and suggestions. 
He is also grateful for the anonymous referees for 
careful reading and helpful comments to improve the manuscript.

\paragraph*{\bf Data availability}
All data generated during this study are included in this article. We have no conflicts of interest to disclose.


\begin{thebibliography}{99}

\bibitem{ach}
M.J. Ablowitz, S. Chakravarty, H. Hahn, 
Integrable systems and modular forms of level 2. J. Phys. A39 (2006), no. 50, 15341-15353.
\bibitem{bot}
J.S. Balakrishnan, K. Ono, W.L. Tsai, Even values of Ramanujan's tau-function. Matematica 1 (2022), no. 2,  395-403.
\bibitem{bgps}
M.A. Bennett, A. Gherga, V. Patel, S. Siksek, 
Odd values of the Ramanujan tau function. Math. Ann. 382 (2022), no. 1-2, 203-238.
\bibitem{br}
S. Brundaban, B. Ramakrishnan, Identities for the Ramanujan tau function and certain convolution
 sum identities for the divisor functions. Ramanujan Math. Soc. Lect. Notes Ser., 23, Ramanujan Math. Soc., Mysore, 2016, 63-75.
\bibitem{chanchua}
H.H. Chan, K.S. Chua, 
Representations of integers as sums of 32 squares. 
Ramanujan J., Rankin memorial issues, vol. 7 (2003),  79-89.
\bibitem{de}
P. Deligne, La conjecture de Weil I. Inst. Hautes Etudes Sci. Publ. Math. (1974), no. 43, 273-307.
\bibitem{hahn}
H. Hahn, convolution sums of some functions on divisors. Rocky Mt. J. Math., Vol. 37, {5}  (2007), 1593-1622.
\bibitem{hahn2}
H. Hahn, Eisenstein series associated with $\Gamma_0(2)$. Ramanujan J., 15 (2008), no. 2, 235-257.
\bibitem{hosw}
J.G. Huard, Z.M. Ou, B.K. Spearman, K.S. Williams, Elementary evaluation of certain convolution sums involving divisor functions in {Number theory for the Millennium}, Vol. II, 2002, 229-274.
\bibitem{ik}
O. Imamo$\bar{\text{g}}$lu, W. Kohnen, Representations of integers as sums of an even number of squares. 
Math. Ann. 333 (2005), no. 4, 815-829.
\bibitem{kk}
M. Kaneko, M. Koike, On modular forms arising from a differential equation of hypergeometric type. Ramanujan J. 7 (2003), no. 1-3, 145-164.
\bibitem{le}
D.H. Lehmer, The vanishing of Ramanujan's $\ta(n)$. Duke Math. J. 1 (1947), 429-433.
\bibitem{ll}
K. Lakein, A. Larsen, Some remarks on small values of $\tau(n)$. Arch. Math. 117 (2021), no. 6, 635-645.
\bibitem{maier}
R.S. Maier, Nonlinear differential equations satisfied by certain classical modular forms. Manuscripta Math., 
vol. 134 (2011), no. 1-2, 1-42.
\bibitem{milne}
S.C. Milne, Hankel determinants of Eisenstein series,
Symbolic computation, number theory, special functions, physics and combinatorics, Gainesville, FL, 
1999, Dev. Math., 4, Kluwer Acad. Publ., Dordrecht, 2001, 171-188.
\bibitem{milne2}
S.C. Milne, New infinite families of exact sums of squares formulas, Jacobi elliptic functions, and Ramanujan's tau function. 
Proc. Nat. Acad. Sci. U.S.A. 93 (1996), no. 26, 15004-15008.
\bibitem{mordell}
L.J. Mordell, On Mr. Ramanujan's empirical expansions of modular functions. Proc. Camb. Phil. Soc. 19 (1917), 117-124.
\bibitem{mor}
P. Moree, On some claims in Ramanujan's `unpublished' manuscript on the partition and tau functions. Ramanujan J. 8 (2004), no. 3, 317-330.
\bibitem{mms}
M.R. Murty, V. K. Murty, T. N. Shorey, Odd values of the Ramanujan  $\tau$-function. 
Bull. Soc. Math. France 115 (1987), no. 3, 391-395.
\bibitem{nikdelan}
Y. Nikdelan, Ramanujan-type systems of nonlinear ODEs for $\Gam_{0}(2)$ and $\Gam_{0}(3)$. Exp. Math. 40 (2022), 409-431.
\bibitem{ramanujan}
S. Ramanujan, On certain arithmetical functions. {Trans. Camb. Phil. Soc.} {22} (1916), 159-184.
\bibitem{shen}
L.C. Shen, On the logarithmic derivative of a theta function and a fundamental identity of Ramanujan. J. Math. Anal. Appl. 177 (1993), no. 1, 299-307.
\bibitem{toh}
P.C. Toh,
Differential equations satisfied by Eisenstein series of level 2. Ramanujan J., 25 (2011), no. 2, 179-194.
\end{thebibliography}
\end{document}